\documentclass{article}

\usepackage{amsmath}
\usepackage{amssymb}
\usepackage{amsthm}

\usepackage{hyperref}
\usepackage{faktor}
\usepackage{tikz-cd}
\usepackage{relsize}

\newtheorem{thm}{Theorem}[section]

\newtheorem*{nthm*}{Theorem 1.5}
\newtheorem*{nnthm*}{Theorem 1.6}

\newtheorem{lem}[thm]{Lemma}
\newtheorem{prop}[thm]{Proposition}
\newtheorem{cor}[thm]{Corollary}

\newtheorem{defn}[thm]{Definition}

\title{The Density of Tuples Restricted by Relatively $r$-Prime Conditions}
\author{Brian D. Sittinger, Vickie V. Chen}
\begin{document}
\maketitle

\begin{abstract}
In order to consider $j$-wise relative $r$-primality conditions that do not necessarily require all $j$-tuples of elements in a Dedekind domain to be relatively $r$-prime, we define the notion of $j$-wise relative $r$-primality with respect to a fixed $j$-uniform hypergraph $H$. This allows us to provide further generalisations to several results on natural densities not only for a ring of algebraic integers $\mathcal{O}$, but also for the ring $\mathbb{F}_q[x]$.
\end{abstract}

\normalsize

\section{Introduction}

In 1976, Benkoski proved that the natural density of the set of relatively $r$-prime $m$-tuples of positive integers (with $rm > 1$) equals $1/\zeta(rm)$, where $\zeta$ is the Riemann zeta function \cite{Benkoski}. This acted as a culmination of the work of Mertens \cite{Mertens}, Lehmer \cite{Lehmer}, and Gegenbauer \cite{Gegenbauer}. Thereafter, T\'oth \cite{Toth} and Hu \cite{Hu1} found the natural density of the set of $j$-wise relatively prime $m$-tuples of positive integers (where $j \leq m$). Extensions of these results have been made to ideals in a ring of algebraic integers $\mathcal{O}$ by Sittinger \cite{Sittinger1, Sittinger3} and subsequently to elements in a ring of algebraic integers as well by Micheli \cite{Ferraguti} and Sittinger \cite{Sittinger2}. Moreover, Morrison \cite{Morrison} and Guo \cite{Guo} gave analogous results for elements in $\mathbb{F}_q[x]$.

\vspace{.1 in}

We can further generalise the notion of $j$-wise relatively primality by considering relative primality conditions that require some but not all $j$-tuples to be relatively prime. A first step in this direction was investigated by Hu \cite{Hu2}, who used graphs to notate which pairs of integers are to be relatively prime.

\begin{defn}
Let $G$ be a simple undirected graph whose $m$ vertices are the positive integers $a_1, \dots, a_m$. We say that the $m$-tuple $(a_1, \dots, a_m)$ is $G$\textbf{-wise} \textbf{relatively prime} if $\gcd(a_i, a_j)=1$ for all adjacent vertices $a_i$ and $a_j$.
\end{defn}

In the case that $G = K_m$, a complete graph with $m$ vertices, we see that $m$ positive integers being $G$-wise relatively prime is the same as saying that these integers are pairwise relatively prime. If we let $i_k(G)$ denote the number of independent sets of $k$ vertices in $G$ (such a set has no two vertices adjacent in $G$), then Hu \cite{Hu2} proved that the density of the set of $G$-wise relatively prime ordered $m$-tuples of positive integers equals
$$\prod_{p} \Big[\sum_{k=0}^{m} i_k(G) \Big(1 - \frac{1}{p}\Big)^{m-k} \Big(\frac{1}{p}\Big)^k\Big],$$
\noindent where the product is over all prime numbers.

\vspace{.1 in}

We now extend the definition of $G$-wise relative primality to accommodate relative primality conditions on certain sets of $j$-tuples of elements from $m$ elements (where $2 \leq j \leq m$) not only from the ring of integers, but also from any Dedekind domain. First, we give a notion of $m$-tuples of elements from a Dedekind domain being relatively $r$-prime.

\begin{defn}
Let $D$ be a Dedekind domain. Fix $r,m \in \mathbb{N}$. We say that $\beta_1, ..., \beta_m \in D$ are \textbf{relatively $r$-prime} if $\mathfrak{p}^r \nmid \langle \beta_1, ..., \beta_m \rangle$ for any prime ideal $\mathfrak{p} \subseteq D$.
\end{defn}

In order to properly generalise the notion of $G$-wise relative primality, we use the concept of a $j$-uniform hypergraph $H$, in which any edge connects exactly $j$ vertices.

\begin{defn}
Let $D$ be a Dedekind domain. Fix $r, j, m \in \mathbb{N}$ where $j \leq m$, and let $H$ be a simple undirected $j$-uniform hypergraph whose $m$ vertices are $\beta_1, \dots, \beta_m \in D$. We say that $\beta_1, \dots, \beta_m \in D$ are $H$\textbf{-wise} \textbf{relatively $r$-prime} if any $j$ adjacent vertices of $H$ are relatively $r$-prime.
\end{defn}

A few remarks are now in order. First, although we state the definitions in this generality, we are in particular interested in the cases of a ring of algebraic integers as well the polynomial rings $\mathbb{F}_q[x]$. Next, suppose we take $D = \mathbb{Z}$, $j = 2$, and $r=1$. Then our hypergraph is a graph $G$, and Definition 1.3 reduces to $m$ integers are $G$-wise relatively prime as defined in \cite{Hu2}. Moreover when $D = \mathcal{O}$ and $H = K^{(j)}_m$, the complete $j$-uniform hypergraph on $m$ vertices, this definition reduces to $m$ elements being $j$-wise relatively $r$-prime as defined in \cite{Sittinger2}.

\begin{defn}
Given a $j$-uniform hypergraph $H$, we say that a subset $S$ of vertices from $H$ is an \textbf{independent vertex set} if $S$ does not contain any hyperedge of $H$. Moreover for any non-negative integer $k$, we let $i_k(H)$ denote the number of independent sets of $k$ vertices in $H$.
\end{defn}

We now state the main results of this article, starting with the algebraic integer case.

\begin{thm}
Fix $r, j, m \in \mathbb{N}$ such that $j \leq m$ and $rm \geq 2$, and let $K$ be an algebraic number field over $\mathbb{Q}$ with ring of integers $\mathcal{O}$.
Then, the density of the set of $H$-wise relatively $r$-prime ordered $m$-tuples of elements in $\mathcal{O}$ equals
$$\prod_{\mathfrak{p}} \Big[\sum_{k=0}^m i_k(H) \Big(1 - \frac{1}{\mathfrak{N}(\mathfrak{p}^r)}\Big)^{m-k} \Big(\frac{1}{\mathfrak{N}(\mathfrak{p}^r)}\Big)^k\Big],$$
\noindent where the product is over all nonzero prime ideals in $\mathcal{O}$.
\end{thm}

\noindent After setting up the pertinent notation in Section 2, we prove Theorem 1.5.

Since the arithmetic in the rings $\mathbb{Z}$ and $\mathbb{F}_q[x]$ have striking similarities (for further details, see \cite{Rosen}), we would expect that we can derive a $H$-wise relatively $r$-prime density statement for $\mathbb{F}_q[x]$. In Section 3, we state and prove an analogue of Theorem 1.4 for the function field case $\mathbb{F}_q[x]$.  

\begin{thm}
Fix $r, j, m \in \mathbb{N}$ such that $j \leq m$ and $rm \geq 2$. Then the density of the set of $H$-wise relatively $r$-prime ordered $m$-tuples of polynomials in $\mathbb{F}_q[x]$ equals
$$\prod_{f \text{ irred.}} \left[\sum_{k=0}^m i_k(H) \Big(1-\dfrac{1}{q^{r\deg f}} \Big)^{m-k}  \Big(\dfrac{1}{q^{r\deg f}} \Big)^k \right],$$
where it is understood that the product is over all monic irreducible polynomials in $\mathbb{F}_q[x]$.
\end{thm}

\noindent \textbf{Remark:} By noting that $\mathfrak{N}(f) =|\mathbb{F}_q[x]/\langle f \rangle| = q^{\deg{f}}$, the analogy between this latter density statement and the one given in the algebraic number ring case is made clear. 

\section{Density of \texorpdfstring{$H$-wise relatively $r$-prime elements in $\mathcal{O}$}{TEXT}}

Let $K$ be an algebraic number field of degree $n$ over $\mathbb{Q}$ with $\mathcal{O}$ as its ring of integers having integral basis $\mathcal{B} = \{\alpha_1, ..., \alpha_n\}$. As a way to generalise the notion of all positive integers less than or equal to some positive constant $M$, we define
$$\mathcal{O}_{\mathcal{B}}[M] = \Big\{\sum_{i=1}^n c_i \alpha_i \, : \, c_i \in [-M, M) \cap \mathbb{Z}\Big\}.$$

The goal of this section is to derive a $H$-wise relatively prime density statement in $\mathcal{O}$ by using the methods developed by \cite{Ferraguti} and \cite{Sittinger2}. First, we define a notion of density for a subset $T$ of $\mathcal{O}^m$ that reduces to the classic notion of density over $\mathbb{Z}$ as follows.

\begin{defn} Let $T \subseteq \mathcal{O}^m$ and fix an integral basis $\mathcal{B}$ of $\mathcal{O}$. The \textbf{upper and lower densities of $T$ with respect to $\mathcal{B}$} are respectively defined as
$$\overline{\mathbb{D}}_{\mathcal{B}}(T) = \limsup_{M \to \infty} \frac{|T \cap \mathcal{O}_{\mathcal{B}}[M]^m|}{|\mathcal{O}_{\mathcal{B}}[M]^m|} \text{ and }\underline{\mathbb{D}}_{\mathcal{B}}(T) = \liminf_{M \to \infty} \frac{|T \cap \mathcal{O}_{\mathcal{B}}[M]^m|}{|\mathcal{O}_{\mathcal{B}}[M]^m|}.$$

\noindent If $\overline{\mathbb{D}}_{\mathcal{B}}(T) = \underline{\mathbb{D}}_{\mathcal{B}}(T)$, we say that its common value is called the \textbf{density of $T$ with respect to $\mathcal{B}$} and denote this as $\mathbb{D}_{\mathcal{B}}(T)$. Whenever this density is independent of the chosen integral basis $\mathcal{B}$, we denote this density as $\mathbb{D}(T)$.
\end{defn}

Although the manner in which we cover $\mathcal{O}$ could potentially depend on the choice of the given integral basis $\mathcal{B}$, it is a direct corollary to Theorem 1.5 that the density of the set of $H$-wise relatively $r$-prime elements in $\mathcal{O}$ is actually independent of the integral basis used.

\vspace{.1 in}

For the remainder of this section, let $S$ be a finite set of rational primes, and fix positive integers $r, j, m$ such that $j \leq m$. Fix a $j$-uniform hypergraph $H$, and define $E_S$ to be the set of $m$-tuples $z = (z_1, \dots, z_m)$ in $\mathcal{O}^m$ such that any ideal generated by $j$ entries of $z$ is $H$-wise relatively $r$-prime with respect to all $\mathfrak{p} \mid \langle p \rangle$ for each $p\in S$. That is, $E_S$ consists of the $H$-wise relatively $r$-prime $m$-tuples of algebraic integers from $\mathcal{O}$ with respect to $S$.

\vspace{.1 in}

In order to aide us in analysing $E_S$, let
$$\pi: \mathcal{O}^m \to \Big(\prod_{\substack{\mathfrak{p} \mid \langle p \rangle \\ p \in S}} \mathcal{O}/\mathfrak{p}^r\Big)^m$$
be the surjective homomorphism induced by the family of natural projections
$${\mathlarger \pi}_{\mathfrak{p}^r}: \mathcal{O} \to \mathcal{O}/\mathfrak{p}^r \text{ for all $\mathfrak{p} \mid \langle p \rangle$ where $p \in S$}.$$

\noindent From the definition of $H$-wise relative $r$-primality of algebraic integers, we immediately deduce the following lemma.

\begin{lem}
For a given prime ideal $\mathfrak{p} \mid \langle p \rangle$ where $p\in S$ and $k \in \{1, 2, \dots, m\}$, let $A_k^{(\mathfrak{p})}$ denote the set of elements in  $(\mathcal{O}/\mathfrak{p}^r)^m$ where exactly $k$ of their $m$ components are 0, and these $k$ components form an independent vertex set in $H$. Then,
$$E_S = \pi^{-1} \Big(\prod_{\substack{\mathfrak{p} \mid \langle p \rangle \\ p \in S}} \bigcup_{k=0}^{m} A_{k}^{(\mathfrak{p})} \Big).$$
\end{lem}

\begin{prop}
Suppose that $\mathfrak{p}$ is a prime ideal in $\mathcal{O}$ that lies above a fixed rational prime $p$, and let $D_p=\sum_{\mathfrak{p} \mid \langle p \rangle}f_p $ where $f_p$ denotes the inertial degree of $\mathfrak{p}$.  If we fix $q \in \mathbb{N}$ and set $N = \prod_{p \in S} p^r$, then
$$\left| E_{S} \cap\mathcal{O}_{\mathcal{B}} [qN]^m \right| = (2q)^{mn}\,\prod_{\substack{\mathfrak{p} \mid \langle p \rangle \\ p \in S}} p^{rm(n-D_p)} \Big[\sum_{k=0}^m i_k(H) \big( \mathfrak{N}(\mathfrak{p}^r) - 1\big)^{m-k} \mathfrak{N}(\mathfrak{p}^r)^k \Big].$$
\end{prop}

\noindent \textbf{Proof:} We first examine the map $\pi$. For brevity, we set $R_p = \prod_{\mathfrak{p} \mid \langle p \rangle} \mathcal{O}/\mathfrak{p}^r$. Then we let $\pi_{\scriptscriptstyle N}$ denote the reduction modulo $N$ homomorphism, and $\psi = (\psi_p)_{p \in S}$ where $\psi_p : (\mathcal{O}/\langle p \rangle^r)^m \to R_p^m$ is the homomorphism induced by the projection maps $\mathcal{O}/\langle p \rangle^r \to R_p$. Finally, let $\overline{\psi}$ be its extension to $(\mathcal{O}/\langle N \rangle)^m$ (by applying the Chinese Remainder Theorem to the primes in $S$). These maps are related to each other through the following diagram 

$$\begin{tikzcd}
 \mathcal{O}^m \arrow{r}{\pi_{\scriptscriptstyle N}} & (\mathcal{O}/\langle N \rangle)^m \arrow{r}{\overline{\psi}} \arrow[swap]{d}{\cong} & (\prod_{p \in S} R_p)^m \arrow{d}{=} \\
& (\prod_{p \in S} \mathcal{O}/\langle p^r \rangle)^m \arrow{r}{\psi} & (\prod_{p \in S} R_p)^m
\end{tikzcd}$$

\noindent and it follows that $\pi = \overline{\psi} \circ \pi_{\scriptscriptstyle N}$.

\vspace{.1 in}

To prove this proposition, we start by examining $\psi^{-1}$. Since for each rational prime $p$ the mapping $\psi_p : (\mathcal{O}/\langle p^r \rangle)^m \to R_p^m$ is a surjective free $\mathbb{Z}_{p^r}$-module homomorphism, we have for all $y \in (\prod_{p \in S} R_p)^m$:
$$|\overline{\psi}^{\, -1}(y)| = \prod_{p \in S} |\psi^{\, -1}_p(y_p)| = \prod_{p \in S} |\ker{\psi_p}| = \prod_{p \in S} p^{rm(n-D_p)}.$$

\vspace{.1 in}

Next, we compute $\Big|\pi^{-1}_{\scriptscriptstyle N}(z) \cap \mathcal{O}_{\mathcal{B}}[qN]^m\Big|$.  Given
$\overline{z} = (\overline{z_1}, \dots , \overline{z_m}) \in (\mathcal{O}/\langle N \rangle)^m$, observe that since $\mathcal{O}/\langle N \rangle$ is a free $\mathbb{Z}_{\scriptscriptstyle N}$-module with basis $\{\pi(\alpha_1), \dots, \pi(\alpha_n)\}$, there exist unique $c_t^j \in [0, N) \cap \mathbb{Z}$ such that
$$\overline{z_j} = \sum_{t=1}^n c_t^j \pi(\alpha_t),$$
\noindent Then for $z = (z_1, \dots, z_m) \in \mathcal{O}^m$, it follows that $\pi_{\scriptscriptstyle N}(z) = \overline{z}$ if and only if
$$z_j = \sum_{t=1}^n (c_t^j + l_t^j N) \alpha_t$$
\noindent for some $l_t^j \in \mathbb{Z}$. Moreover, since we need $l_t^j \in [-q, q) \cap \mathbb{Z}$ for each pair of indices $j$ and $t$, we deduce that
$$\Big|\pi^{-1}_{\scriptscriptstyle N}(z) \cap \mathcal{O}_{\mathcal{B}}[qN]^m \Big| = (2q)^{mn}.$$

\vspace{.1 in}

We are ready to compute $\Big|E_S \cap \mathcal{O}_{\mathcal{B}}[qN]^m\Big|$. By the definition of $A_{k}^{(\mathfrak{p})}$, we have for any fixed $k$ and $\mathfrak{p}$:
$$|A_{k}^{(\mathfrak{p})}| = i_k(H) \big(\mathfrak{N}(\mathfrak{p}^r)-1\big)^{m-k}\mathfrak{N}(\mathfrak{p}^r)^k.$$
Since we know from the last lemma that
$E_S = \pi^{-1}(J)$, where
$$J = \psi^{-1} \Big(\prod_{\substack{\mathfrak{p} \mid \langle p \rangle \\ p \in S}} \bigcup_{k=0}^m A_{k}^{(\mathfrak{p})} \Big),$$
it immediately follows that
$$|J| =\prod_{\substack{\mathfrak{p} \mid \langle p \rangle \\ p \in S}} p^{rm(n-D_p)}  \sum_{k=0}^{m} i_k(H) \big(\mathfrak{N}(\mathfrak{p}^r)-1\big)^{m-k}\mathfrak{N}(\mathfrak{p}^r)^k.$$

\noindent Therefore, we conclude that
\begin{align*}
\Big|E_S \cap \mathcal{O}_{\mathcal{B}}[qN]^m\Big|
&=(2q)^{mn} |J|\\
&=(2q)^{mn}\,\prod_{\substack{\mathfrak{p} \mid \langle p \rangle \\ p \in S}} p^{rm(n - D_p)}
\Big[\sum_{k=0}^m i_k(H) \big(\mathfrak{N}(\mathfrak{p}^r)-1\big)^{m-k}\mathfrak{N}(\mathfrak{p}^r)^k \Big], \end{align*}	

\noindent as desired. $\blacksquare$

\vspace{.1 in}

\noindent We now compute the density of $E_S$.

\begin{lem}
Using the previous notation, we have for any integral basis $\mathcal{B}$ of $\mathcal{O}$,
$$\mathbb{D}(E_S) = \mathbb{D}_{\mathcal{B}}(E_S) = \prod_{\substack{\mathfrak{p} \mid \langle p \rangle \\ p \in S}}
\Big[\sum_{k = 0}^m i_k(H) \Big(1 - \frac{1}{\mathfrak{N}(\mathfrak{p}^r)}\Big)^{m-k} \Big(\frac{1}{\mathfrak{N}(\mathfrak{p}^r)}\Big)^k\Big].$$
\end{lem}

\noindent \textbf{Proof:} Define the sequence $\{a_j\}$ by
$\displaystyle a_j = \frac{|E_S \cap \mathcal{O}_{\mathcal{B}}[j]^m|}{|\mathcal{O}_{\mathcal{B}}[j]^m|}$, and let $D$ denote value of the density in question.

\vspace{.1 in}

First, we consider the subsequence $\{a_{qN}\}_{q \in \mathbb{N}}$, where $N = \prod_{p \in S} p^r$. We claim that this subsequence is constant.
By the previous proposition along with the definitions for $N$ and $D_p$,
\begin{align*}
a_{qN}	& = \frac{1}{(2qN)^{mn}} \Big[(2q)^{mn} \cdot \prod_{\substack{\mathfrak{p} \mid \langle p \rangle \\ p \in S}} p^{rm(n - D_p)} \sum_{k = 0}^m i_k(H) \Big(\mathfrak{N}(\mathfrak{p}^r) - 1\Big)^{m-k} \mathfrak{N}(\mathfrak{p}^r)^k\Big]\\
	& = \prod_{\substack{\mathfrak{p} \mid \langle p \rangle \\ p \in S}} \Big[\sum_{k = 0}^m i_k(H) \Big(1 - \frac{1}{\mathfrak{N}(\mathfrak{p}^r)}\Big)^{m-k} \Big(\frac{1}{\mathfrak{N}(\mathfrak{p}^r)}\Big)^k\Big].
\end{align*}	

\noindent Hence, $\{a_{qN}\}$ is a constant subsequence and converges to $D$.

\vspace{.1 in}

Next, we show that $\{a_{c+qN}\}$ also converges to $D$ for any $c \in \{1, 2, \dots, N-1\}$, we first find bounds for $a_{c+qN}$. To this end, note that



$$a_{qN} \Big(\frac{2qN}{2c + 2qN}\Big)^{mn} \leq a_{c+qN} \leq a_{(q+1)N} \Big(\frac{2(q+1)N}{2c + 2qN}\Big)^{mn}.$$

\noindent By letting $q \to \infty$ and applying the Squeeze Theorem, we conclude that $\{a_{c+qN}\}$ converges to $D$ for any $c \in \{1, 2, \dots, N-1\}$. Finally, since $\{a_{c+qN}\}$ converges to $D$ for any $c \in \{0, 1, \dots, N-1\}$, we conclude that $\{a_j\}$ converges to $D$. $\blacksquare$

\vspace{.1 in}

Note that the density in Lemma 2.4 is independent of the integral basis $\mathcal{B}$ used. Now we are ready to establish to the main theorem of this section. For convenience, we restate it here before proving it.

\begin{nthm*}
Fix $r, j, m \in \mathbb{N}$ such that $j \leq m$ and $rm \geq 2$, and let $K$ be an algebraic number field over $\mathbb{Q}$ with ring of integers $\mathcal{O}$.
Then, the density of the set $E$ consisting of $H$-wise relatively $r$-prime ordered $m$-tuples of elements in $\mathcal{O}$ equals
$$\prod_{\mathfrak{p}} \Big[\sum_{k=0}^m i_k(H) \Big(1 - \frac{1}{\mathfrak{N}(\mathfrak{p}^r)}\Big)^{m-k} \Big(\frac{1}{\mathfrak{N}(\mathfrak{p}^r)}\Big)^k\Big],$$
\noindent where the product is over all nonzero prime ideals in $\mathcal{O}$.
\end{nthm*}

\noindent \textbf{Proof:} Fix $t \in \mathbb{N}$ and let $S_t$ denote the set of the first $t$ rational primes. For brevity, we write $E_t = E_{S_t}$. Since $E_t \supseteq E$,
$$\overline{\mathbb{D}}_{{\mathcal{B}}}(E) \leq \overline{\mathbb{D}}_{{\mathcal{B}}}(E_t) = \mathbb{D}(E).$$
\noindent Observe that the last equality is due to the existence of $\mathbb{D}(E)$. Letting $t \to \infty$, $$\mathbb{D}_{{\mathcal{B}}}(E) \leq \prod_{\mathfrak{p}} \Big[\sum_{k = 0}^m i_k(H) \Big(1 - \frac{1}{\mathfrak{N}(\mathfrak{p}^r)}\Big)^{m-k} \Big(\frac{1}{\mathfrak{N}(\mathfrak{p}^r)}\Big)^k \Big].$$

It remains to show the opposite inequality. Noting that
$\mathbb{D}_{\mathcal{B}}(E_t) - \overline{\mathbb{D}}_{\mathcal{B}}(E_t \backslash E) \leq \underline{\mathbb{D}}_{\mathcal{B}}(E)$, it suffices to show that $\lim_{t \to \infty} \overline{\mathbb{D}}_{\mathcal{B}}(E_t \backslash E) = 0$.

\vspace{.1 in}

To this end, we introduce the following notation.
Let $\mathfrak{p}$ be a prime ideal in $\mathcal{O}$, $p_t$ be the $t^{\text{th}}$ rational prime, and $M$ be a positive integer.

\begin{itemize}

\item[(1)] We write $\mathfrak{p} \succ M$ iff $\mathfrak{p}$ lies over a rational prime greater than $M$.

\item[(2)] We write $M \succ \mathfrak{p}$ iff the rational prime lying under $\mathfrak{p}$ is less than $M$.

\end{itemize}

\noindent Using this notation, we can write
$$E_t \backslash E \subseteq \bigcup_{\mathfrak{p} \succ p_t} \Big(\prod_{j=1}^m \mathfrak{p}^r \Big) \subseteq \mathcal{O}^m,$$
where it is understood that $\prod_{j=1}^m \mathfrak{p}^r$ is the subset of $\mathcal{O}^m$ such that each entry of the $m$-tuple is an element of $\mathfrak{p}^r$. Then, we see that
$$(E_t \backslash E) \cap \mathcal{O}_{\mathcal{B}}[M]^m \subseteq \bigcup_{CM^n \succ \mathfrak{p} \succ p_t} \; \prod_{j=1}^m \Big(\mathfrak{p}^r  \cap \mathcal{O}_{\mathcal{B}}[M]\Big)$$

\noindent for some constant $C > 0$ dependent only on $\mathcal{B}$, and thus

$$\overline{\mathbb{D}}_{\mathcal{B}} (E_t \backslash E) \leq \limsup_{M \to \infty}  \sum_{CM^n \succ \mathfrak{p} \succ p_t} |(\mathfrak{p}^r  \cap \mathcal{O}_{\mathcal{B}}[M])^m| \cdot (2M)^{-mn}.$$

\noindent By \cite[Proposition 13]{Ferraguti}, there exist constants $c, d > 0$ independent of $M$ and $\mathfrak{p}$ such that
$$|(\mathfrak{p}^r  \cap \mathcal{O}_{\mathcal{B}}[M])^m| \leq \frac{(2M)^{mn}}{\mathfrak{N}(\mathfrak{p}^r)^m} + c \Big(\frac{2M}{d \, \mathfrak{N}(\mathfrak{p}^r)^{1/n}} + 1 \Big)^{mn-1}.$$

\noindent Using this bound along with the facts that $\mathfrak{N}(\mathfrak{p}) \geq p$ for every $\mathfrak{p}$ lying above a fixed rational prime $p$, and at most $n$ prime ideals lie above a fixed rational prime, we obtain

\begin{align*}
\overline{\mathbb{D}}_{\mathcal{B}} (E_t \backslash E) & \leq \limsup_{M \to \infty}  \sum_{CM^n \succ \mathfrak{p} \succ p_t}  \Big[ \frac{1}{\mathfrak{N}(\mathfrak{p}^r)^m} + c \Big(\frac{2M}{d \, \mathfrak{N}(\mathfrak{p}^r)^{1/n}} + 1 \Big)^{mn-1} (2M)^{-mn}\Big]\\
& \leq \limsup_{M \to \infty}  \sum_{CM^n > p > p_t}  \Big[ \frac{n}{p^{rm}} + cn \Big(\frac{2M}{d p^{r/n}} + 1 \Big)^{mn-1}  (2M)^{-mn}\Big].
\end{align*}

\vspace{.1 in}

It remains to show that the right side goes to 0 as $t \to \infty$. First, observe that for all sufficiently large $M$, we have $2M/{d p^{r/n}} > 1$ and thus
$$\Big(\frac{2M}{d p^{r/n}} + 1 \Big)^{mn-1} (2M)^{-mn} < \Big(\frac{2}{d}\Big)^{mn} \cdot \frac{1}{p^{rm}}.$$

\noindent Then, by writing $\displaystyle A = n + cn(2/d)^{mn}$ which is a constant independent of $M$ and $p$, we deduce that
$$\overline{\mathbb{D}}_{\mathcal{B}} (E_t \backslash E) \leq  \limsup_{M \to \infty} \sum_{CM^n > p > p_t} \frac{A}{p^{rm}} \leq \sum_{k=p_t}^{\infty} \frac{A}{k^{rm}}$$
for all sufficiently large $M$.

\noindent Finally since $\displaystyle\sum_{k=1}^{\infty} \frac{1}{k^{rm}}$ is convergent, we conclude that $\overline{\mathbb{D}}_{\mathcal{B}} (E_t \backslash E) = 0$. \, $\blacksquare$

\vspace{.1 in}

To conclude this section, we now state a corollary that indicates how this main result provides a generalisation of the work from \cite{Sittinger2}.

\begin{cor}
Fix $r, j, m \in \mathbb{N}$ such that $j \leq m$ and $rm \geq 2$, and let $K$ be an algebraic number field over $\mathbb{Q}$ with ring of integers $\mathcal{O}$. Then the density of the set of $j$-wise relatively $r$-prime ordered $m$-tuples of elements in $\mathcal{O}$ equals
$$\prod_{\mathfrak{p}} \Big[\sum_{k=0}^{j-1} \binom{m}{k} \Big(1 - \frac{1}{\mathfrak{N}(\mathfrak{p}^r)}\Big)^{m-k}  \Big(\frac{1}{\mathfrak{N}(\mathfrak{p}^r)}\Big)^k\Big].$$
\end{cor}

\noindent \textbf{Proof:} Take $H = K_m^{(j)}$ as the hypergraph, and observe that 
$$i_k(H) = \begin{cases} \binom{m}{k} & \text{ if } 0 \leq k \leq j-1\\ 0 & \text{otherwise.} \end{cases}$$ 
Applying Theorem 1.5 immediately yields the desired result. $\blacksquare$

\section{Density of \texorpdfstring{$H$-wise relatively $r$-prime elements in $\mathbb{F}_q[x]$}{TEXT}}

Let $\mathbb{F}_q[x]$ be the ring of polynomials over the finite field $\mathbb{F}_q$ where $q=p^k$ for some prime $p$ and $k\in \mathbb{N}$. The goal of this section is to derive a $H$-wise density statement in $\mathbb{F}_q[x]$ by using methods developed in \cite{Guo}.

\vspace{.1 in}

In order to define a suitable definition of density in $\mathbb{F}_q[x]$, we begin by giving an enumeration of the polynomials in $\mathbb{F}_q[x]$.  Denoting the elements of $\mathbb{F}_q$ as $a_0 = 0$, $a_1, \dots, a_{q-1}$, let $\Sigma$ be the set of all $(a_{d_0}, a_{d_1}, a_{d_2}, \dots)$ whose entries are in $\mathbb{F}_q$ and $d_i = 0$ for all sufficiently large $i$. Then since non-negative integers have a unique expansion base $q$, where $q$ is a positive integer greater than 1, we have a bijection $\Phi :\Sigma \rightarrow \mathbb{Z}_{\geq 0}$ defined by
$$\Phi (a_{d_0}, a_{d_1}, \dots) = \sum_{i=0}^{\infty} d_i q^i.$$

\noindent Using this bijection, we define for each $j\in \mathbb{Z}_{\geq 0}$
$$f_j(x)=\sum_{i=0}^\infty a_{d_i}x^i, \text{ where } j = \phi (a_{d_0}, a_{d_1}, \dots ).$$

\noindent Note that $\mathbb{F}_q[x]=\{ f_j(x) \, : \, j \in \mathbb{Z}_{\geq 0} \}$, thereby giving an ordering of the elements in $\mathbb{F}_q[x]$. Now, we are able to define a density in this ring.

\begin{defn}
Fix a positive integer $m\geq 2$, and let $\mathcal{M}_N$ be the subset of $(\mathbb{F}_q[x])^m$ consisting of $m$-tuples of elements in $\mathbb{F}_q[x]$ whose entries are taken from $\{f_0,f_1, \dots ,f_N \}$.  For any subset $T\subseteq (\mathbb{F}_q[x])^m$, we define the
\textbf{upper and lower densities of $T$} are respectively defined as
$$\overline{\mathbb{D}}(T) = \limsup_{N \to \infty} \frac{|T \cap \mathcal{M}_N|}{|\mathcal{M}_N|} \text{ and }\underline{\mathbb{D}}(T) = \liminf_{N \to \infty} \frac{|T \cap \mathcal{M}_N|}{|\mathcal{M}_N|}.$$

\noindent If $\overline{\mathbb{D}}(T) = \underline{\mathbb{D}}(T)$, we say that its common value is called the \textbf{density of $T$} and denote this as $\mathbb{D}(T)$.
\end{defn}

Let $S$ be a finite set of irreducible polynomials in $\mathbb{F}_q[x]$, and fix $r,j,m \in \mathbb{N}$ satisfying $j \leq m$. Fix a $j$-uniform hypergraph $H$, and let 
$E_S$ denote the set of $m$-tuples of polynomials from $\mathbb{F}_q[x]$ that are $H$-wise relatively $r$-prime with respect to all irreducible polynomials in $S$.

\vspace{.1 in}

For the following lemma and proposition, let
$$\pi: (\mathbb{F}_q[x])^m \to \Big(\prod_{f \in S} \mathbb{F}_q[x]/\langle f^r\rangle\Big)^m$$
\noindent be the surjective homomorphism induced by the family of natural projections
$${\mathlarger \pi}_{f^r}: \mathbb{F}_q[x] \to \mathbb{F}_q[x]/\langle f^r\rangle \text{ for each } f \in S.$$

\noindent As in the algebraic integer case, the following lemma follows immediately from the definition of $H$-wise relative $r$-primality of elements in $\mathbb{F}_q[x]$.

\begin{lem}
For a given irreducible polynomial $f \in S$, let $A_k^{(f)}$ denote the set of elements in  $(\mathbb{F}_q[x]/\langle f^r \rangle)^m$ where exactly $k$ of their $m$ components are 0, and these $k$ components form an independent vertex set in $H$. Then,
$$E_S = \pi^{-1} \Big(\prod_{f \in S} \bigcup_{k=0}^{m} A_{k}^{(f)} \Big).$$
\end{lem}

\begin{prop}
Let $N = bq^{\deg{F}} - 1$ where $b \in \mathbb{N}$, and $F=\prod_{f\in S}f^r$.  Then,
$$ \left| E_S \cap \mathcal{M}_N \right| \, = \,(bq^{\deg{F}})^m \prod_{f \in S} q^{-rm \deg{f}} \cdot \sum_{k = 0}^m i_k(H) (q^{r \deg{f}} - 1)^{m-k} (q^{r \deg{f}})^k.$$
\end{prop}

\noindent \textbf{Proof:} Let $\pi_{F} $ denote the reduction modulo $F$ homomorphism, and let
$$\psi :\, (\mathbb{F}_q[x]/\langle F\rangle )^m \rightarrow \Big(\prod_{f \in S} (\mathbb{F}_q[x]/\langle f^r \rangle \Big)^m \rightarrow \prod_{f \in S} (\mathbb{F}_q[x]/\langle f^r \rangle )^m,$$
\noindent where the first part of $\psi$ is induced by the Chinese Remainder Theorem and the second part is an obvious isomorphism of free $\mathbb{F}_q[x]$-modules.


Now we compute $|\pi^{-1}_{F}\big( h(x) \big) \cap \mathcal{M}_N|$. By the Division Algorithm, we have that
$$\{ f_l (x) \}_{l=0}^N = \{ f_s(x)\cdot x^{\deg{F}} + f_t(x) \mid 0\leq t \leq q^{\deg{F}}-1 \, \textnormal{and} \, 0\leq s \leq b-1  \}.$$
\noindent Then for any fixed $s \in \{0,1, \dots, b-1 \}$, the map $\pi_F$ restricted to
$$\{ f_s(x)\cdot x^{\deg{F}}+f_t(x) \}^{q^{\deg{F}}-1}_{t=0} \rightarrow \mathbb{F}_q[x]/\langle F \rangle $$ is one-to-one. Since $|\ker(\pi_{F})| = b^m$, we conclude that $|\pi^{-1}_{F}\big( h(x) \big) \cap \mathcal{M}_N| = b^m$.

\vspace{.1 in}

We are now ready to compute $|E_S\cap \mathcal{M}_N|$. We know that $E_S = \pi^{-1}(J)$, where
$$\displaystyle J = \psi^{-1} \Big(\prod_{f \in S} \bigcup_{k = 0}^m A_k^{(f)} \Big).$$

\noindent Since for any fixed $k \in \{0, 1, \dots, m\}$ and $f \in S$ we have
$$|A_k^{(f)}| = i_k(H) (q^{r \deg{f}} - 1)^{m-k} (q^{r \deg{f}})^k,$$

\noindent we deduce that
$$|J| = q^{m\deg{F}} \prod_{f \in S} q^{-rm \deg{f}} \cdot \sum_{k = 0}^m i_k(H) (q^{r \deg{f}} - 1)^{m-k} (q^{r \deg{f}})^k.$$

\noindent Therefore,

\begin{align*}
|E_S\cap \mathcal{M}_N| &= b^m \cdot |J|\\
&= (bq^{\deg{F}})^m \prod_{f \in S} q^{-rm \deg{f}} \cdot \sum_{k = 0}^m i_k(H) (q^{r \deg{f}} - 1)^{m-k} (q^{r \deg{f}})^k. \, \blacksquare
\end{align*}

We now find the density of $E_S$.

\begin{lem}
Using the notation from Proposition 3.3,
$$\mathbb{D}(E_S) = \displaystyle\prod_{f\in S} \Big[\sum_{k=0}^m i_{k}(H) \Big(1-\frac{1}{q^{r\deg f}}\Big)^{m-k} \Big(\frac{1}{q^{r\deg f}}\Big)^k \Big].$$
\end{lem}

\noindent \textbf{Proof:} Let $a_j = \displaystyle \frac{|E_S \cap \mathcal{M}_j|}{|\mathcal{M}_j|}$ and let $D$ be the value of the density in question. For notational brevity, we let $n = q^{\deg{F}}$.

\vspace{.1 in}

We first consider the subsequence $\{a_{bn - 1}\}_{b \in \mathbb{N}}$. By Proposition 3.3, we find that
$$\frac{|E_S\cap \mathcal{M}_{bn-1}|}{|\mathcal{M}_{bn-1}|}  =  \prod_{f \in S}\Big[ \sum_{k = 0}^m i_k(H) \Big(1 - \frac{1}{q^{r \deg{f}}} \Big)^{m-k} \Big(\frac{1}{q^{r \deg{f}}}\Big)^k\Big].$$

\noindent Hence, $\{a_{bn - 1}\}$ trivially converges to $D$.

\vspace{.1 in}

Next, we show $\{a_{bn + c}\}$ converges to $D$ as well for each $c \in \{0, 1, \dots, n-2\}$. In a manner reminiscent of the proof to Lemma 2.4, we find that
$$\Big(\frac{bn}{bn + c + 1}\Big)^m a_{bn - 1} \leq a_{bn + c} \leq \Big(\frac{(b+1)n}{(b+1)n + c + 1}\Big)^m a_{(b+1)n - 1}.$$
Letting $b \to \infty$, the Squeeze Theorem implies that $\{a_{bn + c}\}$ converges to $D$ for each
$c \in \{0, 1, \dots, n-2\}$. Finally, since $\{a_{bn + c}\}$ converges to $D$ for each
$c \in \{0, 1, \dots, n-1\}$, we conclude that $\{a_j\}$ converges to $D$, as desired. $\blacksquare$




\vspace{.1 in}

Now we are ready to state and prove the main theorem of this section.

\begin{nnthm*}
Fix $r, j, m \in \mathbb{N}$ such that $j \leq m$ and $rm \geq 2$. Then the density of the set of $H$-wise relatively $r$-prime ordered $m$-tuples of polynomials in $\mathbb{F}_q[x]$ equals
$$\prod_{f \text{ irred.}} \left[\sum_{k=0}^m i_k(H) \Big(1-\dfrac{1}{q^{r\deg f}} \Big)^{m-k}  \Big(\dfrac{1}{q^{r\deg f}} \Big)^k \right],$$
where it is understood that the product is over all monic irreducible polynomials in $\mathbb{F}_q[x]$.
\end{nnthm*}

\noindent \textbf{Proof:} Fix a monic irreducible polynomial $f \in \mathbb{F}_q[x]$ and let $K_f$ denote the set of ordered $m$-tuples $(g_1, \dots, g_m)$ such that $f$ divides the gcd of $k$ of the entries from $(g_1, \dots, g_m)$ whenever these $k$ entries form an independent vertex set. Then by Lemma 3.4, we have
$$\mathbb{D}(K_f) = 1-\sum_{k=0}^m i_k(H) \Big( 1-\dfrac{1}{q^{r\deg f}} \Big)^{m-k} \Big(\dfrac{1}{q^{r\deg f} } \Big)^k.$$

\noindent However for any $x\in[0,1]$, Bernoulli's Inequality implies that
\begin{align*}
\sum_{k=0}^m i_k(H) x^k (1-x)^{m-k} &\geq (1-x)^m + mx(1-x)^{m-1}\\
&= (1-x)^{m-1} (1 + (m-1)x)\\
& \geq (1 - (m-1)x)(1 + (m-1)x)\\
&= 1 - (m-1)^2 x^2.
\end{align*}

\noindent Therefore, letting $x = q^{-\deg{f}}$ yields
$$\mathbb{D}(K_f) \leq \Big(\dfrac{m-1}{q^{r\deg f}}\Big)^2.$$

Next, let $S_t$ be the set of monic irreducible polynomials of a degree greater or equal to $t$ where $t\in\mathbb{N}$, and set $E_t=E_{S_{t}}$.
Moreover, let $\hat{S}$ be the set of all monic irreducible polynomials in $\mathbb{F}_q[x]$. Then,
\begin{align*}
\overline{\mathbb{D}}(E_t \backslash E) &\leq \limsup_{N\to \infty} \dfrac{|(\bigcup_{f\in{\hat{S}\backslash S_t}}K_f)\cap \mathcal{M}_N|}{|\mathcal{M}_N|}\\
&\leq \limsup_{N\to \infty} \dfrac{\sum_{f\in \hat{S}\backslash S_t}|K_f\cap \mathcal{M}_N|}{|\mathcal{M}_N|}\\
&\leq \sum_{f\in \hat{S}\backslash S_t} \overline{\mathbb{D}}(K_f).
\end{align*}

\noindent Since $\overline{\mathbb{D}}(K_f) = \mathbb{D}(K_f)$, we obtain
\begin{align*}
\overline{\mathbb{D}}(E_t \backslash E)
&\leq \sum_{f\in \hat{S}\backslash S_t}\mathbb{D}(K_f)\\
&\leq \sum_{f\in \hat{S}\backslash S_t} \Big( \dfrac{m-1}{q^{r\deg f} } \Big)^2\\
&=\sum_{j=t+1}^{\infty} \dfrac{(m-1)^2}{q^{2rj}}\cdot \varphi(j),
\end{align*}
\noindent where $\varphi(j)$ denotes the number of monic irreducible polynomials of degree $j$ in $\mathbb{F}_q[x]$.

\vspace{.1 in}

Since any irreducible polynomial over $\mathbb{F}_q[x]$ with degree $j$ divides $x^{q^j}-x$ (which has no multiple roots), we have $j\cdot \varphi(j)\leq q^j$.  Therefore
$$\overline{\mathbb{D}}(E_t \backslash E) \leq \sum_{j=t+1}^{\infty} \dfrac{(m-1)^2}{jq^{(2r-1)j}} \leq \dfrac{(m-1)^2}{q^t(q-1)},$$

\noindent in which the last inequality follows from
\begin{align*}
\sum_{j=t+1}^{\infty}\dfrac{1}{jq^{(2r-1)j}}
&=\dfrac{1}{q^{(2r-1)(t+1)}}\cdot \sum_{j=0}^{\infty}\dfrac{1}{(j+t+1)q^{(2r-1)j}}\\
&\leq \dfrac{1}{q^{(2r-1)(t+1)}}\cdot \sum_{j=0}^{\infty}\dfrac{1}{q^{(2r-1)j}}\\
&\leq \dfrac{1}{q^t (q-1)}.
\end{align*}

Next, since $E\cap \mathcal{M}_N\subseteq E_t\cap \mathcal{M}_N$, it follows that
$$\overline{\mathbb{D}}(E) \leq \overline{\mathbb{D}}(E_t) \leq \mathbb{D}(E_t).$$

\noindent Similarly, since $E\cap \mathcal{M}_N=(E_t \cap \mathcal{M}_N)-\big((E_t\backslash E) \cap \mathcal{M}_N\big)$, we obtain
\begin{align*}
\underline{\mathbb{D}}(E) &\geq \underline{\mathbb{D}}(E) - \overline{\mathbb{D}}(E\backslash E_t)\\
&\geq \mathbb{D}(E_t)-\dfrac{(m-1)^2}{q^t(q-1)}.
\end{align*}

\noindent Finally noting that $\mathbb{D}(E_t)$ exists, we conclude by letting $t\rightarrow \infty$ that
\begin{align*}
\mathbb{D}(E) &=\lim_{t\to\infty}\mathbb{D}(E_t)\\
&= \lim_{t\to\infty}\prod_{f\in S_t}\Big[\sum_{k=0}^m i_k(H) \Big(1-\dfrac{1}{q^{r\deg f}} \Big)^{m-k} \Big(\dfrac{1}{q^{r\deg f}} \Big)^k\Big]\\
&=\prod_{f \text{ irred.}} \Big[\sum_{k=0}^m i_k(H) \Big(1-\dfrac{1}{q^{r\deg f}} \Big)^{m-k} \Big(\dfrac{1}{q^{r\deg f}} \Big)^k\Big],
\end{align*}
and this concludes the proof. $\; \blacksquare$

\vspace{.1 in}

In a manner reminiscent of the previous section, we conclude by giving without proof the analogue of Corollary 2.5 for $\mathbb{F}_q[x]$ as originally given in \cite{Guo}. 

\begin{cor}
Fix $r, j, m \in \mathbb{N}$ such that $j \leq m$ and $rm \geq 2$. Then the density of the set of $j$-wise relatively $r$-prime ordered $m$-tuples of elements in $\mathbb{F}_q[x]$ equals
$$\prod_{f \, \text{irred.}} \Big[\sum_{k=0}^{j-1} \binom{m}{k} \Big(1 - \frac{1}{q^{r\deg f}}\Big)^{m-k}  \Big(\frac{1}{q^{r\deg f}}\Big)^k\Big],$$
where it is understood that the product is over all monic irreducible polynomials in $\mathbb{F}_q[x]$.
\end{cor}

\end{document}